\documentclass[a4paper]{article}
\usepackage{amsmath}

\newcommand{\p}{\partial}
\newcommand{\cD}{{\cal D}}
\newcommand{\cA}{{\cal A}}

\newcommand{\cS}{{\cal S}}
\newcommand{\cC}{{\cal C}}
\newcommand{\cF}{{\cal F}}
\newcommand{\cP}{\cal P}
\newcommand{\cU}{{\cal U}}
\newcommand{\cE}{\cal E}
\newcommand{\cT}{{\cal T}}
\newcommand{\cX}{\cal X}

\newcommand{\lp}{\left(}
\newcommand{\rp}{\right)}
\newcommand{\raa}{\rightarrow}
\newcommand{\ap}{\alpha}
\newcommand{\ga}{\gamma}

\newcommand{\De}{\Delta}

\newcommand{\Ci}{C^{\infty}}
\newcommand{\E}{\ell}
\newcommand{\m}{\!\!\mid}
\newcommand{\N}{{\mathbf N}}

\newcommand{\R}{{\mathbf R}}
\newcommand{\op}[1]{\!\!\mathop{\rm ~#1}\nolimits}

\mathchardef\za="710B  
\mathchardef\zb="710C  
\mathchardef\zg="710D  
\mathchardef\zd="710E  
\mathchardef\zve="710F 
\mathchardef\zz="7110  
\mathchardef\zh="7111  
\mathchardef\zvy="7112 
\mathchardef\zi="7113  
\mathchardef\zk="7114  
\mathchardef\zl="7115  
\mathchardef\zm="7116  
\mathchardef\zn="7117  
\mathchardef\zx="7118  
\mathchardef\zp="7119  
\mathchardef\zr="711A  
\mathchardef\zs="711B  
\mathchardef\zt="711C  
\mathchardef\zu="711D  
\mathchardef\zvf="711E 
\mathchardef\zq="711F  
\mathchardef\zc="7120  
\mathchardef\zw="7121  
\mathchardef\ze="7122  
\mathchardef\zy="7123  
\mathchardef\zf="7124  
\mathchardef\zvr="7125 
\mathchardef\zvs="7126 
\mathchardef\zf="7127  
\mathchardef\zG="7000  
\mathchardef\zD="7001  
\mathchardef\zY="7002  
\mathchardef\zL="7003  
\mathchardef\zX="7004  
\mathchardef\zP="7005  
\mathchardef\zS="7006  
\mathchardef\zU="7007  
\mathchardef\zF="7008  
\mathchardef\zW="700A  

\newcommand{\be}{\begin{equation}}
\newcommand{\ee}{\end{equation}}
\newcommand{\rar}{\rightarrow}

\newcommand{\bea}{\begin{eqnarray}}
\newcommand{\eea}{\end{eqnarray}}
\newcommand{\beas}{\begin{eqnarray*}}
\newcommand{\eeas}{\end{eqnarray*}}

\begin{document}

\title{Derivations of the\\ Lie algebras of differential operators
\footnote{This work was supported by MCESR Grant RD/C.U.L./02-010
and by KBN, grant No 2 P03A 020 24.}}
\author{J. Grabowski, N. Poncin}\maketitle

\newtheorem{re}{Remark}
\newtheorem{theo}{Theorem}
\newtheorem{prop}{Proposition}
\newtheorem{lem}{Lemma}
\newtheorem{cor}{Corollary}
\newtheorem{ex}{Example}

\begin{abstract}
This paper encloses a complete and explicit description of the
derivations of the Lie algebra $\cD(M)$ of all linear differential
operators of a smooth manifold $M$, of its Lie subalgebra
$\cD^1(M)$ of all linear first-order differential operators of
$M$, and of the Poisson algebra $\cS(M)=\op{Pol}(T^*M)$ of all
polynomial functions on $T^*M,$ the symbols of the operators in
$\cD(M).$ It turns out that, in terms of the Chevalley cohomology,
$H^1(\cD(M),\cD(M))=H^1_{\op{DR}}(M)$,
$H^1(\cD^1(M),\cD^1(M))=H^1_{\op{DR}}(M)\oplus\R^2$, and
$H^1(\cS(M),\cS(M))=H^1_{\op{DR}}(M)\oplus\R$. The problem of
distinguishing those derivations that generate one-parameter
groups of automorphisms and describing these one-parameter groups
is also solved.
\end{abstract}

\section{Introduction}

In \cite{PS}, Pursell and Shanks proved the well-known result stating that
the Lie algebra of all smooth compactly supported vector fields of a
smooth manifold characterizes the differentiable structure of the variety.
Similar upshots were obtained in numerous subsequent papers dealing with
different Lie algebras of vector fields and related algebras (see e.g.
\cite{A,Am,AG,JG,JG3,HM,O,S}).

Derivations of certain infinite-dimensional Lie algebras arising
in Geometry were also studied in different situations (note that
in infinite dimension there is no such a clear correspondence
between derivations and one-parameter groups of automorphisms as
in the finite-dimensional case). Let us mention a result of
L.~S.~Wollenberg \cite{Wol} who described all derivations of the
Lie algebra of polynomial functions on the canonical symplectic
space $\R^2$ with respect to the Poisson bracket. It turned out
that there are outer derivations of this algebra in contrast to
the corresponding Weyl algebra. This can be viewed as a variant of
a "no-go" theorem (see \cite{J}) stating that the Dirac
quantization problem \cite{Dir} cannot be solved satisfactorily
because the classical and the corresponding quantum algebras are
not isomorphic as Lie algebras. An algebraic generalization of the
latter fact, known as the {\it algebraic "no-go" theorem,} has
been proved in \cite{GG} by different methods. Derivations of the
Poisson bracket of all smooth functions on a symplectic manifold
have been determined in \cite{ADML} (for the real-analytic case,
see \cite{Gr2}). Another important result is the one by F.~Takens
\cite{Tak} stating that all derivations of the Lie algebra
$\cX(M)$ of smooth vector fields on a manifold $M$ are inner. The
same turned out to be valid for analytic cases \cite{Gr1}. Some
cases of the Lie algebras of vector fields associated with
different geometric structures were studied in a series of papers
by Y.~Kanie \cite{Ka1}--\cite{Ka4}.

Our work \cite{GP} contains Shanks-Pursell type results for the
Lie algebra $\cD(M)$ of all linear differential operators of a
smooth manifold $M$, for its Lie subalgebra $\cD^1(M)$ of all
linear first-order differential operators of $M$, and for the
Poisson algebra $\cS(M)=\op{Pol}(T^*M)$ of all polynomial
functions on $T^*M,$ the symbols of the operators in $\cD(M).$
Furthermore, we computed all the automorphisms of these algebras
and showed that $\cD(M)$ and $\cS(M)$ are not integrable. The
current paper contains a description of their derivations, so it
is a natural continuation of this previous work and can be
considered as a generalization of the results of Wollenberg and
Takens. It is also shown which derivations generate one-parameter
groups of automorphisms and the explicit form of such
one-parameter groups is provided.

\section{Notations and definitions}

Throughout this paper, $M$ is as usually assumed to be a smooth,
Hausdorff, second countable, connected manifold of dimension $n$.

Recall that the space $\cD(M)$ (or $\cD$ for short) of linear
differential operators on $\Ci(M)$ (or $\cA$ for short) is
filtered by the order of differentiation, $\cD^{i}$ being the
space of at most $i$-th order operators (for $i\ge 0$;
$\cD^{i}=\{0\}$ for $i<0$), and is equipped with an associative
and so a Lie algebra structure, $\circ$ and $[.,.]$ respectively,
such that $\cD^{i}\circ\cD^j\subset \cD^{i+j}$ and
$[\cD^{i},\cD^j]\subset\cD^{i+j-1}$. Obviously, $\cD^0=\cA$ is an
associative commutative subalgebra and $\cD^1$ is a Lie subalgebra
of $\cD.$ 
We denote by $\cD_c$ (respectively $\cD^{i}_c$) the algebra of
differential operators (respectively the space of (at most) $i$-th
order operators) that vanish on constants. For instance $\cD^1_c$
is the Lie algebra ${\cX}(M)=\op{Vect}(M)$ (or $\cX$ for short) of
vector fields of $M$, i.e. the Lie algebra $\op{Der}\cA$ of
derivations of the algebra of functions. Observe also that we have
the canonical splittings $\cD=\cA\oplus\cD_c$,
$\cD^{i}=\cA\oplus\cD^{i}_c.$

The classical counterpart of $\cD$, the space $\cS(M)$ (or $\cS$
for short) of symmetric contravariant tensor fields on $M$, is of
course naturally graded, $\cS_i$ being the space of $i$-tensor
fields (for $i\ge 0$; $\cS_i=\{0\}$ for $i<0$). 
This counterpart $\cS$ is isomorphic--- even as a
$\cD^1_c$-module---to the space $\op{Pol}(T^*M)$ of smooth
functions on $T^*M$ that are polynomial on the fibers.
Furthermore, it is a commutative associative and a Poisson
algebra. These structures $\cdot$ and $\{.,.\}$ verify
$\cS_i\cdot\cS_j\subset\cS_{i+j}$ and
$\{\cS_i,\cS_j\}\subset\cS_{i+j-1}$ respectively. The Poisson
bracket can be viewed as the symmetric Schouten bracket or the
standard symplectic bracket. Note that $\cS_0=\cA$ is an
associative and Lie-commutative subalgebra of $\cS$. Clearly,
$\cS$ is filtered by $\cS^{i}=\oplus_{j\le i}\cS_j$ and $\cS^1$ is
a Lie subalgebra of $\cS$ isomorphic to $\cD^1$ and $\cA\oplus
{\cal X}$.

The algebras $\cD$ and $\cS$ are models of a quantum and a classical
Poisson algebra in the sense of \cite{GP}. All the results of this paper
apply to these algebras. It is well known that $\cD^{i},$ $\cS_i,$ and
$\cS^{i}$ are algebraically characterized in the following way:
\begin{equation}\{D\in\cD:[D,\cA]\subset\cD^{i}\}=\cD^{i+1}\;(i\ge
-1)\label{characdf},\end{equation}
\begin{equation}\{S\in\cS:\{S,\cA\}\subset\cS_i\}=\cA+\cS_{i+1}\;(i\ge
-1), \label{characsfg}\end{equation} and
\begin{equation}\{S\in\cS:\{S,\cA\}\subset\cS^{i}\}=\cS^{i+1}\;(i\ge
-1)\label{characsf}.\end{equation} Moreover, $\cS$ is the
classical algebra induced by the quantum algebra $\cD$. Thus,
$\cS_i=\cD^{i}\slash\cD^{i-1}$. For any non-zero $D\in\cD$, the
degree $\op{deg}(D)$ of $D$ is the lowest $i$, such that
$D\in\cD^{i}\backslash\cD^{i-1}$. If $\op{cl}_j$ is the class in
the quotient $\cS_j$, the (principal) symbol $\sigma(D)$ of $D$ is
\[\sigma(D)=\op{cl}_{\op{deg}(D)}(D)\] and the symbol $\sigma_i(D)$ of order $i\ge \op{deg}(D)$
is defined by \[\sigma_i(D)=\op{cl}_i(D)=\begin{cases} 0,\mbox{ if }i>\op{deg}(D), &\\
\sigma(D),\mbox{ if }i=\op{deg}(D).&\end{cases}\] Then, the
commutative multiplication and the Poisson bracket of $\cS$ verify
\begin{equation}\sigma(D_1)\cdot\sigma(D_2)=\sigma_{\op{deg}(D_1)+\op{deg}(D_2)}(D_1\circ
D_2)\;\;(D_1,D_2\in\cD)\label{point-circ}\end{equation} and
\begin{equation}\{\sigma(D_1),\sigma(D_2)\}=\sigma_{\op{deg}(D_1)+
\op{deg}(D_2)-1}([D_1,D_2])\;\;(D_1,D_2\in\cD).\label{Poisson-Lie}\end{equation}

\section{Locality and weight}

The characterizations (\ref{characdf}), (\ref{characsfg}), and
(\ref{characsf}) of the filters $\cD^{i+1}$ of $\cD$ and the terms
${\cS}_{i+1}$ and filters ${\cS}^{i+1}$ of $\cS$ ($i\ge -1$), can
be "extended" in the following way:

\begin{lem} For any $i\ge -1$ and any $k\ge 1$, we have
\begin{equation}\{D\in\cD:[D,\cD^k]\subset\cD^{i}\}=\R\cdot 1+
\cD^{i-k+1},\label{characddk}\end{equation}
\begin{equation}\{S\in\cS:\{S,\cS_k\}\subset\cS_i\}=\R\cdot 1+
\cS_{i-k+1},\label{characsskg}\end{equation} and
\begin{equation}\{S\in\cS:\{S,\cS_k\}\subset\cS^{i}\}=\R\cdot 1+
\cS^{i-k+1}.\label{characssk}\end{equation}\label{charack}
\end{lem}

\textit{Proof.} (i) Note first that $\{S\in\cS:\{S,\cS_k\}=0\}=\R\cdot 1$.
Of course, we need only show that the commutation of $S$ with $\cS_k$
implies $S\in\R\cdot 1.$ But this is obvious: on a connected Darboux chart
domain $U$, take for instance the polynomials $S_k\in\cS_k$ defined by
$S_k(x;\xi)=(\xi_i)^k$ and $S_k(x;\xi)=x^i(\xi_i)^k$ ($x\in U,$ $\xi\in
(\R^n)^*$, $\{\xi_j,x^i\}=\delta^i_j$, $i,j\in\{1,\ldots,n\}$).

More generally, we have $\{S\in\cS:\{S,\cS_k\}\subset\cS_i\}=\R\cdot
1+\cS_{i-k+1}$ for all $i\ge -1$. Take $i\ge 0$. Writing $S=S_{i-k+1}+S'$
with $S_{i-k+1}\in\cS_{i-k+1}$ and $S'\in\cS\ominus\cS_{i-k+1}$, we get
$\{S',\cS_k\}\subset\cS_i\cap\lp\cS\ominus\cS_i\rp$, so $\{S',\cS_k\}=0$
and $S'\in\R\cdot 1$. Hence the conclusion.

(ii) In order to prove (\ref{characssk}), observe that it is enough to
consider the case $i\ge 0.$ If $\{S,\cS_k\}\subset\cS^{i}$ and
$S=\sum_jS_j$, $S_j\in\cS_j$, we have $\{S_j,\cS_k\}=0$ for all $j>i-k+1$.
So $S_j\in\R\cdot 1$ and $S\in\R\cdot 1+\cS^{i-k+1}$.

(iii) Assume $[D,\cD^k]\subset\cD^{i}$ and
$D\in\cD\backslash\cD^{i-k+1}$, so that $\op{deg}(D)>i-k+1$.
Clearly, $\sigma_k:\cD^k\raa\cS_k$ is surjective, so any
$S_k\in\cS_k\backslash \{0\}$ reads $S_k=\sigma(\Delta)$,
$\op{deg}(\Delta)=k$. Thus,
\[\{\sigma(D),S_k\}=\{\sigma(D),\sigma(\Delta)\}=\sigma_{\op{deg}(D)+k-1}([D,\Delta])=0.\]
So $\sigma(D)\in\R\cdot 1$ and $D\in\R\cdot 1$. Eventually, $D\in\R\cdot
1+\cD^{i-k+1}$.
\rule{1.5mm}{2.5mm}\\

Let $(\cP,[.,.])$ be either the Lie algebra $(\cD,[.,.])$, its Lie
subalgebra $(\cD^1,[.,.])$, or the Poisson algebra
$(\cS,\{.,.\}).$ The sign "$\cdot$" stands for the multiplication
"$\!\circ$" of differential operators and the multiplication
"$\cdot$" of polynomials of $T^*M$. We denote by $\op{Der}\cP$ the
Lie algebra of all derivations of the Lie algebra $(\cP,[.,.])$.

\begin{prop}Any derivation of the Lie algebra $\cP$ is a
local operator.\end{prop}

\textit{Proof.} If $P\in {\cal P}^{i}$ vanishes on an open
$U\subset M$ and if $x_0\in U$, we have \[P=\sum_k[X_k,P_k]\] for
certain $X_k\in {\cal X}$, $P_k\in{\cal P}^{i}$ with
$X_k\m_V=P_k\m_V=0$ for some neighborhood $V\subset U$ of $x_0$.
In the quantum case, this follows for instance from \cite{NP}. In
the classical case, a straightforward adaptation of \cite[Ex.
$12$]{DWLc1} shows that the set $\{L_X:\Gamma(\cS^{\le
i}TM)\raa\Gamma(\cS^{\le i}TM)\;|\;X\in {\cal X}\},$ where
$\Gamma(\cS^{\le i}TM)$ is the space of smooth sections of the
tensor bundle $\oplus_{j\le i}\cS^jTM,$ is locally transitive.
Since it is obviously stable under locally finite sums, the
announced result is a direct consequence of \cite[Prop. $3$, Def.
$2$]{DWLc1}.

For any $C\in \op{Der}\cP$, the preceding decomposition of $P$ and
the derivation property then imply that $(CP)(x_0)=0$.
\rule{1.5mm}{2.5mm}

\begin{lem} There is a finite set
$\cF=\{f_1,\ldots,f_m\}\subset\Ci(M)$
($m\le 2n+1$), such that \\
(j) the $\Ci(M)$-module $\Omega^1(M)$ of differential $1$-forms on
$M$ is spanned by $d\cF=\{df_1,\ldots,df_m\}$,\\
(jj) if $P\in\cP$ verifies $[P,\cF]\subset {\cal P}^{i}$ then
$P\in {\cal P}^{i+1}$, for any $i\ge -1$.\label{charac}
\end{lem}


\textit{Proof.} Assertion (j) is a consequence of Whitney's embedding
theorem (see \cite[Prop. $2.6$]{AG}, \cite{Wh}). It suffices to prove (jj)
for $i=-1$. Indeed, by induction, if (jj) is verified for $i$ ($i\ge -1$)
then it is for $i+1$:
\[\begin{array}{ll}[P,\cF]\subset {\cal P}^{i+1}&\Rightarrow [[P,\cF],\cA]
\subset {\cal P}^{i}\\
&\Rightarrow [[P,\cA],\cF]\subset {\cal P}^{i}\\
&\Rightarrow [P,\cA]\subset {\cal P}^{i+1}\\
&\Rightarrow P\in {\cal P}^{i+2}.\end{array}\] As
\[[D,\cF]=0\Rightarrow \{\sigma(D),\cF\}=0, \forall D\in\cD,\] it is enough to consider
the classical case, which is obvious in view of (j). Indeed, if
$f\in\cA$ we have $df=\sum_{s=1}^mg_s df_s$ ($g_s\in\cA$) and if
$S\in\cS=\cS(M)=\op{Pol}(T^*M)$ and $\Lambda$ denotes the
canonical Poisson tensor of $T^*M$, then
\[\{S,f\}=\Lambda(dS,df)=\sum_{s=1}^mg_s\{S,f_s\}=0.\] Hence the result. \rule{1.5mm}{2.5mm}

\begin{prop} Any derivation $C$ of the Lie algebra $\cP$ has a bounded weight,
i.e. there is $d\in\N,$ such that
\[C({\cal P}^{i})\subset {\cal P}^{i+d},\forall i\in\N.\]\end{prop}
\textit{Proof.} Set
$d=\op{max}\{\op{deg}(Cf_s),s\in\{1,\ldots,m\}\}$, where
$\op{deg}$ is the degree in the filtered algebra $\cP$ and where
the set $\cF=\{f_1,\ldots,f_m\}$ is that of Lemma \ref{charac}.
Then $C$ maps all functions into ${\cal P}^d$. Indeed, if
$f\in\cA$ we have for any $f_s\in\cF$,
\[0=C[f,f_s]=[Cf,f_s]+[f,Cf_s]\] and $[Cf,f_s]\in {\cal P}^{d-1},$ so that
$Cf\in {\cal P}^d.$ The announced result can then once more be obtained by
induction. Take $P\in {\cal P}^{i+1}$ ($i\ge 0$) and apply again the
derivation property:
\[C[P,f_s]=[CP,f_s]+[P,Cf_s],\forall s\in\{1,\ldots,m\}.\] Hence
the conclusion. \rule{1.5mm}{2.5mm}\\

\textbf{Remark}: Evidently, for $C\in \op{Der}\cD^1$, we have
$C(\cD^0)\subset\cD^1$ and $C(\cD^1)\subset\cD^1\subset\cD^2$.

\section{Corrections by inner derivations}

\begin{prop} Let $C\in \op{Der}\cP$. There is (a non-unique) $P\in\cP,$ such that
$C-\op{ad}\,P\in \op{Der}\cP$ respects the filtration. The set of
all elements of $\cP$ that have this property is then $P+{\cal
P}^1$. \label{filt}\end{prop}

\textit{Proof.} Take an arbitrary derivation $C$ of the Lie algebra $\cP$.
Let $(U_{\iota},\varphi_{\iota})_{\iota\in I}$ be an atlas of $M$ and
${\cal U}=(U,\varphi=(x^1,\ldots,x^n))$ any chart of this atlas. The
restriction $C\m _U$ of the local operator $C$ to the domain $U$ is of
course a derivation of the Lie algebra ${\cal P}_U$, similar to $\cP$ but
defined on $U$ instead of $M$.

Set now
\[P_C^{{\cal U},i}=C\m _U(x^{i})\in {\cal P}_U^d.\] This element
$P_C^{{\cal U},i}$ is equal to or symbolically represented by a polynomial
of $T^*U$ of type
\[P_C^{{\cal U},i}(x;\xi)=\sum_{\mid\ap\mid\le
d}\ga_{\ap}^{i}(x)\,\xi^{\ap},\] where we used standard notations,
$\ga_{\ap}^{i}\in \Ci(U)$ and $\xi\in(\R^n)^*.$

Since it follows from $C\m _U[x^{i},x^j]=0$ that $[P_C^{{\cal
U},i},x^j]=[P_C^{{\cal U},j},x^{i}],$ we get
\[\p_{\xi_j}P_C^{{\cal U},i}(x;\xi)=\p_{\xi_i}P_C^{{\cal U},j}(x;\xi).\] Thus, there is
a polynomial of $T^*U$,
\[P_C^{\cU}(x;\xi)=\sum_{\mid\ap\mid\le
d+1}\ga_{\ap}(x)\xi^{\ap}\] (polynomial character in $\xi$ and smooth
dependence on $x$ easily checked), such that
\[\p_{\xi_i}P_C^{\cU}(x;\xi)=P_C^{{\cal U},i}(x;\xi),\forall i\in\{1,\ldots,n\}.\]
Finally, $P_C^{\cU}\in {\cal P}^{d+1}_U$ (interpret---if necessary---the
polynomial as differential operator) and
\[C\m _U(x^{i})=[P_C^{\cU},x^{i}],\forall i.\]
For any function $f\in\cA$ and any $i\in\{1,\ldots,n\}$, we then have
\[\begin{array}{ll}
0=C\m_U[f\m_U,x^{i}]&=[(Cf)\m_U,x^{i}]+[f\m_U,[P_C^{\cU},x^{i}]]\\
&=[(Cf)\m_U-[P_C^{\cU},f\m_U],x^{i}].\end{array}\] In view of Lemma
\ref{charac}, this entails that
\begin{equation}(Cf)\m_U-[P_C^{\cU},f\m_U]\in\Ci(U).\label{locfunc}\end{equation}

Now we will glue together the elements $P_C^{\cU}\in {\cal P}^{d+1}_U$.
Let $(U_{\iota},\varphi_{\iota},\psi_{\iota})_{\iota\in I}$ be a partition
of unity subordinated to the considered atlas and set
\[P_C=\sum_{\iota}\psi_{\iota}P_C^{{\cal U}_{\iota}}.\]
Clearly $P_C\in{\cal P}^{c+1}.$ Furthermore, $C-\op{ad}\,P_C$ is a
derivation of $\cP=\cD$ and of $\cP=\cS$. Let us emphasize that
for ${\cal P}=\cD^1$, this map $C-\op{ad}\,P_C$ verifies the
derivation property in $\cD^1$, but is a priori only linear from
$\cD^1$ into $\cD$. For any $\cP$, it respects the filtration.
Indeed, for any $f\in\cA$ and any open $V\subset M$, we have
\begin{equation}\begin{array}{ll}\lp Cf-[P_C,f]\rp\m_V&=(Cf)\m_V-[\sum_i\psi_i\m_VP_C^{{\cal U}_i}\m_{U_i\cap V},f\m_V]\\
&=\sum_i\psi_i\m_V\lp (Cf)\m_{U_i}-[P_C^{{\cal
U}_i},f\m_{U_i}]\rp\m_{U_i\cap
V}\;\;\in\Ci(V)\label{filtf},\end{array}\end{equation} in view of
(\ref{locfunc}). We can now proceed by induction. So assume that
$CP-[P_C,P]\in{\cal P}^{i},\forall P\in{\cal P}^{i}$ ($i\ge 0$). Then, if
$P\in{\cal P}^{i+1}$ and $f\in\cA$,
\begin{equation}\begin{array}{ll}[CP-[P_C,P],f]&=C[P,f]-[P,Cf]-[P_C,[P,f]]+[P,[P_C,f]]\\
&=(C[P,f]-[P_C,[P,f]])-[P,Cf-[P_C,f]]\in{\cal
P}^{i}\label{filtd1}.\end{array}\end{equation} Hence the result for
$\cP=\cD$ and $\cP=\cS$. For ${\cal P}=\cD^1$, Equation (\ref{filtd1})
shows that $[P_C,\cD^1]\subset\cD^1$. In view of Lemma \ref{charack}, this
means that $P_C\in\cD^1.$ Finally, Equation (\ref{filtf}) allows to see
that $C(\cD^0)\subset\cD^0$,
so that $C$ respects the filtration.  \rule{1.5mm}{2.5mm}\\

\noindent \textbf{Remark}: Thus the inner derivation of Proposition
\ref{filt} can be taken equal to $0$, in the case ${\cal P}=\cD^1.$

\begin{prop} If $C\in \op{Der}\cP$ respects the filtration, there is a unique vector
field $Y\in \op{Der}\cA\subset\cP$ such that $C-\op{ad}\,Y\in
\op{Der}\cP$ respects the filtration and \[\lp
C-\op{ad}\,Y\rp\m_{\cA}=\kappa\;\op{id},\] where $\kappa\in\R$ is
uniquely determined by $C$.\label{restrfunct}
\end{prop}

\textit{Proof.} Consider a derivation $C$ of $\cP$ that respects
the filtration and denote by $\cC(\cP)$ the centralizer of
$\op{ad}\,\cA$ in the Lie algebra ${\cal E}=\op{End}\,\cP$ of
endomorphisms of $\cP$, i.e. the Lie subalgebra $\cC({\cal
P})=\{\psi\in{\cal E}:[\psi,\op{ad}\,{\cal A}]_{\cE}=0\},$ where
$[.,.]_{\cE}$ is the commutator of endomorphisms of $\cP.$ The
derivation $C\in \op{Der}\cP$ induces a derivation $\op{ad}_{\cal
E}C\in \op{Der}\cal E$, which respects the centralizer. Indeed,
for any $\psi\in\cC(\cP)$, we have
\[\begin{array}{ll}[\lp \op{ad}_{\cal E}C\rp(\psi),\op{ad}\,\cA]_{\cal
E}&=[[C,\psi]_{\cal E},\op{ad}\,\cA]_{\cal
E}\\&=-[[\psi,\op{ad}\,\cA]_{\cal E},C]_{\cal
E}-[\psi,[C,\op{ad}\,\cA]_{\cal E}]_{\cal E}=0,\end{array}\] as
$[C,\op{ad}\,f]_{\cal E}=\op{ad}(Cf)\in \op{ad}\cA$ for each
$f\in\cA$, since $C$ is a derivation that respects the filtration.

It follows from the description of the centralizer, see
\cite[Theo. $3$]{GP}, that if $\psi\in\cC(\cP)$, then $\psi$
respects the filtration and there is $\psi_1\in\cC(\cP)$, such
that $\psi_1({\cal P}^{i})\subset{\cal P}^{i-1}$ and
$\psi=\E_{\psi(1)}+\psi_1$. Obviously, the left multiplication
$\E_f:{\cal P}\ni P\raa f\cdot P\in\cP$ by an arbitrary $f\in\cA$
belongs to the centralizer $\cC(\cP)$. So $\lp \op{ad}_{\cal
E}C\rp(\E_f)\in\cC(\cP)$ and for any $g\in\cA$, $[C,\E_f]_{\cal
E}(g)=[C,\E_f]_{\cal E}(1)\cdot g,$ i.e. $\lp C-C(1)\rp (f\cdot
g)=\lp C-C(1)\rp(f)\cdot g+f\cdot\lp C-C(1)\rp(g)$. As constants
are the only central elements and as derivations map central
elements to central elements, $C(1)=\kappa$ ($\kappa\in\R$), and
the preceding result means that $\lp C-\kappa\,\op{id}\rp\m_{\cA}$
is a vector field $Y$. Finally, $C-\op{ad}\,Y\in \op{Der}\cP$
respects the filtration and $\lp
C-\op{ad}\,Y\rp\m_{\cA}=\kappa\,\op{id}\m_{\cA}$. Uniqueness of
$Y$ is readily obtained. Indeed, if $Y$ is a suitable vector
field, then the corresponding constant is necessarily
$\kappa=C(1)$ and $Y$ is unique. \rule{1.5mm}{2.5mm}

\section{Characterization of the derivations for the Lie algebra $\mathbf{\cD^1(M)}$}


Let $\mid\eta\mid$ be a fixed smooth nowhere zero $1$-density. The
associated divergence $\op{div}_{\mid\eta\mid}$ (or simply
$\op{div}$) is defined for any vector field $X$ as the unique
function $\op{div}_{\mid\eta\mid}X$ that verifies
$L_X\!\!\mid\eta\mid=\lp\op{div}_{\mid\eta\mid}X\rp\mid\eta\mid$,
where $L_X\!\!\mid\eta\mid$ is the Lie derivative of
$\mid\eta\mid$ in the direction of $X$. In any local coordinate
system in which $\mid\eta\mid$ is a constant multiple of the
standard density, this divergence reads
$\op{div}_{\mid\eta\mid}X=\sum_i\p_{x^{i}}X^{i}$, with
self-explaining notations. For details regarding the origin of the
class of the divergence, we refer the reader to \cite{Lec}.

\begin{theo} A map $C:\cD^1(M)\raa\cD^1(M)$ is a derivation of the
Lie algebra $\cD^1(M)=\Ci(M)\oplus \op{Vect}(M)$ of first order
differential operators on $\Ci(M)$, if and only if it can be
written in the form
\begin{equation}C_{Y,\zk,\zl,\zw}(X+f)=[Y,X+f]+\kappa\,f+\lambda\,\op{div}X+\omega(X),
\label{formderivd1}
\end{equation}
where $Y\in \op{Vect}(M)$, $\kappa,\lambda\in\R$, and
$\omega\in\Omega^1(M)\cap \op{ker}\,d$. All these objects
$Y,\kappa,\lambda,\zw$ are uniquely determined by $C$.
\label{derivd1}\end{theo}

\begin{cor} The first group of the Chevalley-Eilenberg
cohomology of the Lie algebra $\cD^1(M)$ of first order differential
operators on $\Ci(M)$ with coefficients in the adjoint representation, is
given by
\[H^1(\cD^1(M),\cD^1(M))\simeq \R^2\oplus H^1_{\op{DR}}(M),\]
where $H^1_{\op{DR}}(M)$ stands for the first space of the de Rham
cohomology of M.\label{cohomd1}\end{cor}

\textit{Proof.} Let $C_1$ be a derivation of $\cP$ that respects
the filtration and reduces to $\kappa\,\op{id}$ ($\kappa\in\R$) on
functions. The derivation property, written for $f\in\cA$ and
$X\in {\cal X}$, shows that $C_1({\cal X})\subset\cA$, and written
for $X,Y\in {\cal X}$, it means that $C_1\m_{{\cal X}}$ is a
$1$-cocycle of the Lie algebra ${\cal X}$ canonically represented
upon $\cA$. The cohomology $H({\cX},{\cA})$ is known (see e.g.
\cite{Fuc} or \cite{DWLc2}). Having fixed a divergence on ${\cal
X}$, we get
\[C_1\m_{{\cal X}}=\lambda\,\op{div}+\omega,\]
with $\lambda\in\R$ and $\omega\in\Omega^1(M)\cap \op{ker}\,d$.
Finally,
\[C_1(X+f)=\kappa\,f+\lambda\,\op{div}X+\omega(X),\forall X\in {\cal X},\forall f\in\cA.\]
To prove uniqueness of $Y,\omega,\kappa,\lambda$, it suffices to
write Equation (\ref{formderivd1}) successively for $1\in\cA$,
$f\in\cA$, and $X\in {\cal X}$. Hence Theorem \ref{derivd1}.
Corollary \ref{cohomd1} is clear. Indeed, if $C$ and
$C'=C+\op{ad}\,(Z+h)$ ($Z\in {\cal X}$, $h\in\cA$) are two
cohomologous $1$-cocycles and if we denote by
$(Y,\omega,\kappa,\lambda)$ and $(Y',\omega',\kappa',\lambda')$
the respective unique quadruples, then necessarily
$\kappa'=\kappa,$ $\lambda'=\lambda,$ and $\omega'=\omega-dh$, so
that the map
\[H^1(\cD^1,\cD^1)\ni [C]\raa\kappa+\lambda+[\omega]\in\R^2\oplus
H^1_{\op{DR}}(M)\] is
a well-defined vector space isomorphism.  \rule{1.5mm}{2.5mm}\\

\textbf{Remarks}: $1$. Observe that the preceding proof is valid for the
generic algebra $\cP$, so not only for $\cD^1$ but also
for $\cS$ and $\cD$.\\

$2$. Note that, for $\lambda\ne 0$ and $\omega=df$ ($f\in\cA$),
$\lambda\,\op{div}X+\omega(X)=\lambda\,\op{div}'X$, where
$\op{div}'X=\op{div}X+\lambda^{-1}\omega(X)$ is another divergence.\\

$3$. Denote $C_{Y,0,0,0}=C_Y$, $C_{0,1,0,0}=C_{\cA}$,
$C_{0,0,1,0}=C_{\op{div}}$, $C_{0,0,0,\zw}=C_\zw$. The Lie algebra
structure of $\op{Der}\cD^1$ is determined by the following
commutation relations (the commutators we miss are just 0): \be
[C_Y,C_{Y'}]=C_{[Y,Y']},[C_Y,C_{\op{div}}]=C_{d(\op{div} Y)},
[C_Y,C_{\zw}]=C_{d(\zw(Y))}, [C_{\cA},C_{\op{div}}]=C_{\op{div}},
[C_{\cA},C_\zw]=C_\zw.\ee

\section{Characterization of the derivations for the Lie algebra $\mathbf{\cS(M)}$}

\textbf{Remark}: Let us recall that we mentioned in \cite[Sect.
4]{GP} two specific types of derivations: the canonical derivation
$\op{Deg}\in \op{Der}\cS$, $\op{Deg}:\cS_i\ni S\raa
(i-1)S\in\cS_i$ and the derivation $\overline{\omega}\in
\op{Der}\cP$ implemented by a closed $1$-form $\omega$ of $M.$ If
$U$ is an arbitrary open subset of $M$ and if $\omega\m_U=df_U$,
$f_U\in\Ci(U)$, this cocycle $\overline{\omega}$ is defined by
\[\overline{\omega}(P)\m_U=[P\m_U,f_U],\forall P\in\cP.\] Remark
that, if ${\cal P}=\cD^1$, $\overline{\omega}$ coincides with the
derivation $\omega$. In the case $\cP=\cS$ the lowering derivation
$\overline{\zw}$ can be interpreted as the action of the vertical vector
field $\zw^v$ (the vertical lift of the section $\zw$ of $T^*M$) on
polynomial functions on $T^*M$.

\begin{theo} A map $C:\cS(M)\raa\cS(M)$ is a derivation of the
Lie algebra $\cS(M)$ of all infinitely differentiable functions of $T^*M$
that are polynomial along the fibers, if and only if it is of the form
\begin{equation}C_{P,\zk,\zw}(S)=\{P,S\}+\kappa\,\op{Deg}(S)+\zw^v(S),
\label{formderivs}
\end{equation}
where $P\in\cS(M)$, $\kappa\in\R,$ and $\omega\in\Omega^1(M)\cap
\op{ker}\,d$. Here $\kappa$ is uniquely determined by $C$, but $P$
and $\omega$ are not. The set of all fitting pairs is
$\{(P+h,\omega+dh):h\in\Ci(M)\}$, so we get uniqueness if we
impose that the polynomial function $P$ vanishes on the 0-section
of $T^*M$. \label{derivs}\end{theo}

\begin{cor} The first group of the Chevalley-Eilenberg
cohomology of the Lie algebra $\cS(M)$ of all polynomial functions on
$T^*M$ with coefficients in the adjoint representation, is given by
\[H^1(\cS(M),\cS(M))\simeq \R\oplus H^1_{\op{DR}}(M).\]\label{cohoms}
\end{cor}
\textit{Proof.} If $C_1\in \op{Der}\cS$ respects the filtration
and coincides with $\kappa\,\op{id}$ ($\kappa\in\R$) on $\cA,$ the
proof of Theorem \ref{derivd1} yields
\[C_1(X+f)=\kappa\,f+\lambda\,\op{div}X+\omega(X),\forall X\in
{\cal X},\forall f\in\cA\] ($\lambda\in\R,\omega\in\Omega^1(M)\cap
\op{ker}\,d$). This outcome is apparent since a derivation of
$\cS$ that respects the filtration, restricts to a derivation of
$\cS^1\simeq\cD^1.$ It is easy to check that
$C_2=C_1+\kappa\,\op{Deg}-\overline{\omega}$ has all the
properties of $C_1,$ but verifies in addition
\[C_2(X+f)=\lambda\,\op{div}X,\forall X\in
{\cal X},\forall f\in\cA.\] The derivation property, written for
$S\in\cS$ and $f\in\cA$, shows inductively that $C_2$ is lowering.
It is easily seen that $\zl=0$, i.e. that $C_2\m_{{\cal S}^1}=0$.
Indeed, for any $X\in{\cal X}$ and $f\in {\cal A}$, we have
\begin{equation}\{C_2(X^2),f\}=C_2\{X^2,f\}=2\zl\lp
X^2(f)+X(f)\op{div}X\rp.\end{equation} The left hand side of this
identity is, with respect to $f$, a differential operator of order
$1$, and the right hand side is of order $2$, if $\zl\neq 0.$
Hence $\zl=0.$

Now, we need only check that any derivation $C$ of ${\cal S}$,
which vanishes on ${\cal S}^1\simeq{\cal D}^1$, is identically
zero. Suppose by induction that $C$ vanishes on ${\cal S}^k$,
$k\ge 1.$ We will show that then $C$ also vanishes on ${\cal
S}_{k+1}.$ For $f\in {\cal A},$ $X\in {\cal X},$ and $S\in {\cal
S}_{k+1},$ we have $0=C(\{S,f\})=\{C(S),f\}$ and $C(\{X,S\})=\{X,
C(S)\}$, so $C$ maps ${\cal S}_{k+1}$ into ${\cal A}$ and
intertwines the adjoint action of ${\cal X}$. Hence
$C(\{fX,X^{k+1}\})=fX\lp C(X^{k+1})\rp$ and the map $D_X: {\cal
A}\ni f\raa C(\{fX,X^{k+1}\})\in {\cal A}$ is a differential
operator of order $0$. On the other hand, $D_X(f)=-(k+1)C\lp
X(f)X^{k+1}\rp$ is of order $0$ only if it is just $0$. Thus,
$C\lp X(f)X^{k+1}\rp=0$, for all $f\in {\cal A}$ and all $X\in
{\cal X}$. Let us now work locally. Homogeneous polynomials of
degree $k+1$ on $(\R^n)^*$, $n=\op{dim}M$ are spanned by
$(k+1)$-th powers $X^{k+1}$, $X\in\R^n$, and any function reads
$X(f)$ for any non-vanishing $X\in\R^n$. So polynomials of the
form $X(f)X^{k+1}$ locally span ${\cal S}_{k+1}$ and $C=0.$
This completes the proof of Theorem \ref{derivs}, except for
uniqueness of $\kappa$ and the convenient pairs $(P,\omega)$.
Equation (\ref{formderivs}), written for $S=1\in\cA$, shows that
$\kappa$ necessarily equals $-C(1)$. Setting $S=f\in\cA$, then
$S=X\in {\cal X}$ in this same equation, we get $\lp
C-C(1)\,\op{id}\rp\m_{\cA}=\lp \op{ad}\,P\rp\m_{\cA}$ resp. $\lp
C-\op{ad}\,P\rp\m_{{\cal X}}=\omega.$ So, if $(P',\omega')$ is
another suitable pair we have $\{P'-P,\cA\}=0$, so that $P'=P+h$,
$h\in\cA$. But then,
$\omega'-\omega=\op{ad}\,(P-P')=-\op{ad}\,h=dh$ on ${\cal X}$.
Corollary \ref{cohoms} is now obvious.
\rule{1.5mm}{2.5mm}\\

\textbf{Remark}: Denote $C_{P,0,0}=C_P$, $C_{0,1,0}={\op{Deg}}$,
$C_{0,0,\zw}=\zw^v$. The Lie algebra structure of $\op{Der}\cS$ is
determined by the following commutation relations (the missing
commutators are just 0):
$$
[C_P,C_{P'}]=C_{\{ P,P'\}}\,, [\op{Deg},C_P]=C_{\op{Deg}(P)}\,,
[{\zw}^v,C_P]=C_{\zw^v(P)}\,, [\zw^v,\op{Deg}]=\zw^v.
$$
\section{Characterization of the derivations for the Lie algebra $\mathbf{\cD(M)}$}

\begin{lem} Any derivation $C\in \op{Der}\,\cD(M)$ that respects the filtration
induces a derivation $\tilde{C}\in \op{Der}\,\cS(M)$, which
respects the graduation:
\[\tilde{C}:\cS_i(M)\ni
S\raa\sigma_i(C(\sigma_i^{-1}(S)))\in\cS_i(M),\]for all
$i\in\N$.\label{quantclassderiv}\end{lem}

\begin{theo} A map $C:\cD(M)\raa\cD(M)$ is a derivation of the
Lie algebra $\cD(M)$ of all linear differential operators on $\Ci(M)$ if
and only if it can be written in the form
\begin{equation}C_{P,\zw}(D)=[P,D]+\overline{\omega}(D),
\end{equation}
where $P\in\cD(M)$ and $\omega\in\Omega^1(M)\cap \op{ker}\,d$ are
not unique. Again the appropriate pairs are $(P+h,\omega+dh)$,
$h\in\Ci(M)$, so we get uniqueness if we impose that $P$ is
vanishing on constants. \label{derivd}\end{theo}

\begin{cor} The first cohomology group of the Lie algebra $\cD(M)$ of all linear differential
operators on $\Ci(M)$ with coefficients in the adjoint representation is
isomorphic to the first space of the de Rham cohomology of M:
\[H^1(\cD(M),\cD(M))\simeq H^1_{\op{DR}}(M).\]\label{cohomd}\end{cor}

\textit{Proof.} Lemma \ref{quantclassderiv} is a consequence of the
surjective character of $\sigma_i:\cD^{i}\raa\cS_i$ and Equation
(\ref{Poisson-Lie}) that links the Poisson and the Lie brackets.

Propositions \ref{filt},\ref{restrfunct} and Theorem \ref{derivd1}
show that in order to establish Theorem \ref{derivd}, we can start
with $C_1\in \op{Der}\cD$, such that
$C_1(\cD^{i})\subset\cD^{i},\forall i\in\N$ and
$C_1(X+f)=\kappa\,f+\lambda\,\op{div}X+\omega(X),$ with the usual
notations. When correcting by $\overline{\omega}$, we get a
filtration-respecting derivation $C_2=C_1-\overline{\omega},$
which maps $X+f$ to $C_2(X+f)=\kappa\,f+\lambda\,\op{div}X.$ The
derivation $\tilde{C}_2$ induced on the classical level then
verifies $\tilde{C}_2(X+f)=\kappa\,f.$ Theorem \ref{derivs} now
implies that
\begin{equation}\tilde{C}_2(S)=-\kappa\,\op{Deg}(S),\label{deg}\end{equation}
for all $S\in\cS.$

Let us emphasize that Theorem \ref{derivd} is based upon Theorem
\ref{derivd1} and Theorem \nolinebreak\ref{derivs}, itself built upon
Theorem \ref{derivd1}. The point is that the degree-derivation is not
generated by any canonical quantum derivation. Therefore the proof of
Theorem \nolinebreak\ref{derivd} is a little bit more complicated than
that of Theorem \ref{derivs}.

Observe first that Equation (\ref{deg}) means that, for each $i\in\N$,
\begin{equation}C_2\m_{\cD^{i}}=\kappa (1-i)\,\op{id}+\chi_i,\label{low}\end{equation}
where $\chi_i\in \op{Hom}_{\R}(\cD^{i},\cD^{i-1})$. Indeed, it
entails that, if $D\in\cD^{i}$, $i\in\N$, the operator
$C_2D-\kappa (1-i)D$ has a vanishing $i$-th order symbol. As
easily checked, $\chi_0=0,$
$\chi_1(X+f)=\kappa\,f+\lambda\,\op{div}X$, and
$\chi_if=i\kappa\,f$ ($X\in {\cal X},f\in\cA,i\in\N$). Injecting
now this structure into the derivation property, written for
$D^{i}\in\cD^{i}$ and $\De^j\in\cD^{j}$, $i,j\in\N$, we obtain
\begin{equation}\chi_{i+j-1}[D^{i},\De^j]=[\chi_iD^{i},\De^j]+[D^{i},\chi_j\De^j].\label{derivpropchi}\end{equation}

When using the decomposition $\cD=\cA\oplus\cD_c$, we denote by
$\pi_0$ and $\pi_c$ the projections onto $\cA$ and $\cD_c$
respectively. Furthermore, if $D\in\cD$, we set $D_0=\pi_0D=D(1)$
and $D_c=\pi_cD=D-D(1)$, and if $C\in \op{End}\,\cD$, we set
$C_0=\pi_0\circ C\in \op{Hom}_{\R}(\cD,\cA)$ and $C_c=\pi_c\circ
C\in \op{Hom}_{\R}(\cD,\cD_c)$.

The projections on $\cA$ of Equation (\ref{derivpropchi}), written
for $D_c^{i}\in\cD^{i}_c$ ($i\ge 2$) and $f\in\cA$, then for
$D^{i}_c\in\cD^{i}_c$ and $\De^j_c\in\cD^j_c$ ($i+j\ge 3$), read

\begin{equation}(\chi_{i,c}D^{i}_c)(f)=(i-1)\kappa\,D^{i}_c(f)+\chi_{i-1,0}[D^{i}_c,f]_c\label{derivpropchidf}\end{equation}
and
\begin{equation}\chi_{i+j-1,0}[D^{i}_c,\De^j_c]=D^{i}_c(\chi_{j,0}\De^j_c)-\De^j_c(\chi_{i,0}D^{i}_c)\label{derivpropchidd}\end{equation}
respectively. Exploiting first Equation (\ref{derivpropchidf})
with $i=2$ and $D^{i}_c=Y^2$, $Y\in {\cal X}$, we get the upshot
\begin{equation}\chi_{2,c}Y^2=(2\zl+\kappa)Y^2+2\zl \lp \op{div}Y\rp Y,\label{deprchdf2}\end{equation}
for all $Y\in{\cal X}$. If we apply the second order symbol
$\zs_2$ to both sides of this equation, we see that
\begin{equation}2\lambda+\kappa=0. \label{kl1}\end{equation}

There is an atlas of $M$ such that in each chart
$(U,x^1,\ldots,x^n)$ the divergence takes the classical form,
$\op{div}\lp \sum_iX^{i}\p_{x^{i}}\rp=\sum_i\p_{x^{i}}X^{i}.$ We
work in such a chart and write $\p$ (resp. $f'$ and $B(f),$
 $f\in\Ci(U))$ instead of $\p_{x^1}$ (resp. $\p f$ and
$\chi_{2,0}(f\p^2)$). For $i=1, D^{i}_c=g\p$, $j=2,
\De^j_c=f\p^2$, $f,g\in\Ci(U)$, Equation (\ref{derivpropchidd})
yields
\begin{equation}B(gf'-2fg')=\zl f'g''+g\lp B(f)\rp'.\label{deprchdd12}\end{equation}
In particular, $B(f')=\lp B(f)\rp'$ and $B(g')=0$. But then $B=0$,
$\zl=0$ (see Equation (\ref{deprchdd12})), $\zk=0$ (see Equation
(\ref{kl1})), and $C_2\m_{{\cal D}^1}=0.$

We now proceed by induction and show that $C_2\m_{{\cal
D}^{k+1}}=0$, if $C_2\m_{{\cal D}^k}=0,$ $k\ge 1.$ As
$C_2\m_{{\cal D}^{k+1}}$ only depends on the $(k+1)$-th order
symbol, as $[C_2(D),f]=0,$ and $C_2\lp[X,D]\rp=[X,C_2(D)],$ for
all $D\in {\cal D}^{k+1},$ $X\in {\cal X},$ and $f\in {\cal A}$,
$C_2$ defines a map $\tilde{C}_2:{\cal S}_{k+1}\raa {\cal A}$ that
intertwines the adjoint action of ${\cal X.}$ We have shown in the
proof of Theorem \ref{derivs} that such a map necessarily
vanishes. Hence $C_2=0.$ This completes the proof of Theorem
\ref{derivd}. Indeed, the statement regarding the appropriate
pairs $(P,\omega)$ is obvious.
The same is true for Corollary \ref{cohomd}. \rule{1.5mm}{2.5mm}\\

\textbf{Remarks}: $1.$ Denote $C_{P,0}=C_P$,
$C_{0,\zw}=\overline{\zw}$. The Lie algebra structure of
$\op{Der}\cD$ is determined by the following commutation relations
(the missing commutator is 0):
$$
[C_P,C_{P'}]=C_{
[P,P']}\,,\,[\overline{\zw},C_P]=C_{\overline{\zw}(P)}\ .
$$

$2$. 
Corollary \ref{cohomd1}, Corollary \ref{cohoms}, and Corollary
\ref{cohomd} imply that the first adjoint cohomology spaces of the Lie
algebras $\cD^1,$ $\cS$, and $\cD$ are independent of the smooth structure
of $M$, provided that the topology of $M$ remains unchanged. \\

$3$. It is worth comparing our cohomological results with those
obtained in other recent papers. Let
$D_M=(\op{End}\,\cA)_{\op{loc},c}$ be the Lie algebra of local
endomorphisms of $\cA$ that vanish on constants. A well-known
theorem of Peetre, \cite{JP}, guarantees that these operators are
locally differential. The main theorem of \cite{NP2} asserts that
the first three local cohomology groups $H^p(D_M,\cA)_{\op{loc}}$
($p\in\{1,2,3\}$) of $D_M$ canonically represented upon $\cA$ are
isomorphic to the corresponding groups $H^p_{\op{DR}}(M)$ of the
de Rham cohomology of $M$. In particular,
\[H^1(D_M,\Ci(M))_{\op{loc}}\simeq H^1_{\op{DR}}(M).\] Let us quote from \cite{AAL}
the outcome \[H^1(\op{Vect}(M),{\cal D}(M))\simeq
H^1(\op{Vect}(M),{\cal D}^{i}(M))\simeq\R\oplus
H^1_{\op{DR}}(M),\] for all $i\in\N$.

\section{Integrability of derivations}

In this section we distinguish those derivations that generate
(smooth) one-parameter groups of automorphisms of the Lie algebra
$\cP$  (we will call such derivations \textit{integrable}) and we
find explicit forms of these one-parameter groups of
automorphisms. The smoothness of a curve in $\op{Aut}\cP$ is
defined in the obvious way with relation to the smooth structure
on $M$. For instance, $\zF_t$ is smooth in $\op{Aut}\cD$ if for
any $D\in\cD$ and any $f\in C^\infty(M)$ the induced map
$(t,x)\mapsto\zF_t(D)(f)(x)$ is a smooth function on $\R\times M$,
a curve $\zF_t$ in $\op{Aut}\cS$ is smooth if for any $S\in\cS$
the induced map $(t,y)\mapsto\zF_t(S)(y)$ is a smooth function on
$\R\times T^*M$, etc. In the following all one-parameter groups
will be assumed to be smooth.

Since the group $\op{Diff}(M)$ of smooth diffeomorphisms of $M$ is
embedded in $\op{Aut}\cP$ (see \cite{GP}), a partial problem is
the determination of one-parameter groups of diffeomorphisms.
This, however, is well known and the one-parameter groups of
diffeomorphisms are just flows $\op{Exp}(tY)$ of complete vector
fields $Y$. Note that in general it is hard to decide if a given
diffeomorphism is implemented by a vector field, since
neighbourhoods of identity in the connected component of the group
$\op{Diff}(M)$ are far from being filled up by flows (even in the
case when $M$ is compact and all vector fields are complete (see
\cite{Gr0,Ko,Pal})); that differs $\op{Diff}(M)$ from
finite-dimensional Lie groups.

Before we start the investigation into one-parameter subgroups in
$\op{Aut}\cP$ we have to define the group-analogue of the
divergence, which is important for the case ${\cal P}={\cal
D}^1(M)$. Let us stress that in this paper the divergence is not
an arbitrary $1$-cocycle of vector fields with coefficients into
functions, but a cocycle obtained as described in \cite{GP} from a
nowhere vanishing $1$-density or as depicted in \cite{GMM} from an
odd volume form. These cocycles or divergences form some
privileged cohomology class. We will integrate any such divergence
$\op{div}:{\cX}(M)\rar\Ci(M)$ to a group 1-cocycle
$J:\op{Diff}(M)\rar\Ci(M)$, which is a sort of Jacobian. 
Indeed, if $|\eta|$ is the odd volume form inducing the divergence
and if $\phi\in \op{Diff}(M)$, we have
$\phi^*|\eta|=J(\phi)|\eta|$ for a unique positive smooth function
$J(\phi)$. It is easily verified that if $\phi$ is a
diffeomorphism between two domains of local coordinates and if $f$
and $g$ are the component functions of $|\eta|$ in the
corresponding bases, then locally
\[J(\phi)(x)=\frac{g(\phi(x))}{f(x)}|\op{det}\,\p_x\Phi|,\] where $\Phi$
is the local form of $\phi$. For any $\phi,\psi\in \op{Diff}(M)$,
we clearly have \be\label{J} J(\phi\circ\zc)=\zc^*(J(\phi))\cdot
J(\zc). \ee In particular,
$$J(\phi^{-1})=\frac{1}{J(\phi)\circ\phi^{-1}}\ .
$$
A similar concept may be found under the name of Jacobi
determinant in \cite[Def. 6.5.12]{AMR}. Let us put
$\op{Div}(\phi)=\ln{J(\phi)}.$

\begin{prop} For any $X\in{\cX}(M)$ and $\phi\in \op{Diff}(M)$, we have
\item{(a)} \be\label{div}\phi^*(\op{div}\,\phi_*(X))=\op{div}
X+X(\op{Div}(\phi))\ee and, if $X$ is complete, \item{(b)}
\be\label{Div} \op{Div}(\op{Exp}(tX))=\int_0^t(\op{div} X)\circ
\op{Exp}(sX)ds. \ee
\end{prop}
\textit{Proof.} (a) By definition of the action of $\phi$ on
vector fields and differential forms,
$\phi^*(i_{\phi_*(X)}\vert\eta\vert)=i_X(\phi^*\vert\eta\vert)$,
so that
$$\phi^*(\op{div}_{\vert\eta\vert}\phi_*(X))=\op{div}_{\phi^*\vert\eta\vert}(X).
$$
Since $\phi^*\vert\eta\vert=J(\phi)\vert\eta\vert$, (\ref{div})
follows.

(b) Let us put $F_t=\op{Div}(\op{Exp}(tX))$. It is easy to see
that $F_{t+s}=F_t+F_s\circ \op{Exp}(tX)$, which implies the
differential equation
\begin{equation}\dot{F_t}=X(F_t)+\dot{F_0}.\label{PDE}
\end{equation} Additionally,
we have the initial conditions $F_0=0$ and, due to
$$\dot{F_0}\vert\eta\vert=\frac{d}{dt}|_{t=0}(\op{Exp}(tX))^*\vert\eta\vert=(\op{div}X)\vert\eta\vert,$$
$\dot{F_0}=\op{div}X$.
Applying formally the variation of constant method, we find
$$F_t=(\op{Exp}(tX))^*(\int^t_0\dot{F_0}\circ \op{Exp}(-sX)ds)=
\int_0^t(\op{div}X)\circ \op{Exp}(sX)ds. \mbox{  }$$ It is easily
verified that this integral is really a solution. Equation
(\ref{PDE}) is in fact a PDE of first order, which can be written
in the form
\[L_{\hat{X}}F=\dot{F}_0,\] with $\hat{X}=\p_t-X\in{\cal X}(\R\times M)$.
A well-known consequence of the theorem of Frobenius allows to see
that this equation, completed by the boundary condition $F|_M=0$,
has locally one unique solution. Hence,
\[F_t=\int_0^t(\op{div}X)\circ \op{Exp}(sX)ds.\]

\subsection{The case $\cD^1(M)$}

Theorem $8$ of \cite{GP} states that an endomorphism $\Phi$ of $\cD^1$ is
an automorphism of the Lie algebra $\cD^1$ if and only if it reads
\begin{equation}\Phi_{\phi,K,\zL,\zW}(X+f)=\phi_*(X)+(K\,f+\Lambda\,
\op{div}X+\Omega(X))\circ\phi^{-1},\label{automd1}
\end{equation}
where $\phi$ is a diffeomorphism of $M$, $K,\Lambda$ are
constants, $K\neq 0$, $\Omega$ is a closed $1$-form on $M$, and
$\phi_*$ is the push-forward
\begin{equation}(\phi_*(X))(f)=(X(f\circ\phi))\circ\phi^{-1},\label{phi*}
\end{equation}
all the objects $\phi,K,\Lambda,\Omega$ being uniquely determined by
$\Phi$. The one-parameter group condition
$$\Phi_{\phi_t,K_t,\zL_t,\zW_t}\circ\Phi_{\phi_s,K_s,\zL_s,\zW_s}=
\Phi_{\phi_{t+s},K_{t+s},\zL_{t+s},\zW_{t+s}}$$ gives immediately
$$\phi_{t+s}=\phi_t\circ\phi_s,\quad K_{t+s}=K_t\cdot K_s,\quad
\zL_{t+s}=\zL_t+K_t\cdot\zL_s,
$$
and, in view of (\ref{div}),
\be\label{zW}\zW_{t+s}=K_t\zW_s+\phi_s^*\zW_t+\zL_t\,d(\op{Div}(\phi_s)),
\ee with the initial conditions $\phi_0=\op{id}_M$, $K_0=1$,
$\zL_0=0$, $\zW_0=0$. One solves easily: $\phi_t=\op{Exp}(tY)$,
$K_t=e^{\zk t}$, $\zL_t=\zl\frac{e^{\zk t}-1}{\zk}$ (with
$\frac{e^{\zk t}-1}{\zk}=t$ if $\zk=0$), for some unique complete
vector field $Y$ and some unique real numbers $\zk,\zl$. To solve
(\ref{zW}) we derive the differential equation
$$\dot{\zW}_t=\zk\zW_t+\zl\,d(\op{Div}(\phi_t))+\phi_t^*\zw,
$$
where $\zw=\dot{\zW}_0$. This is an inhomogeneous linear equation,
which can be solved by the method of variation of the constant. We
get
$$\zW_t=e^{\zk t}\int_0^te^{-\zk s}(\zl\,
d(\op{Div}(\phi_s))+\phi^*_s\zw)ds
$$
and, in view of (\ref{Div}),
$$\zW_t=\int_0^te^{\zk(t-s)}\left(\zl\,
d\left(\int_0^s \op{div}Y\circ
\op{Exp}(uY)du\right)+(\op{Exp}(sY))^*\zw\right)ds.
$$
Since $\dot{\phi}_0=Y$, $\dot{K}_0=\zk$, $\dot{\zL}_0=\zl$,
$\dot{\zW}_0=\zw$, we get the following:

\begin{theo} A derivation
\[C_{Y,\zk,\zl,\zw}(X+f)=[Y,X+f]+\kappa\,f+\lambda\,\op{div}X+\omega(X)\] of $\cD^1(M)$ induces
a one-parameter group $\zF_t$ of automorphisms of $\cD^1(M)$ if
and only if the vector field $Y$ is complete. In this case the
group is of the form
\beas&\Phi_t(X+f)=(\op{Exp}(tY))_*(X)+\left(e^{\zk
t}\,f+\,\zl\frac{e^{\zk t}-1}{\zk} \op{div}X\right)\circ
\op{Exp}(-tY)+\cr
&\left(\int_0^te^{\zk(t-s)}\left(\zl\int_0^sX(\op{div}Y\circ
\op{Exp}(uY))du+\lp(\op{Exp}(sY))^*\zw\rp(X)\right)ds\right)\circ
\op{Exp}(-tY).\label{1automd1} \eeas
\end{theo}

\subsection{The case $\cS(M)$}
We know from \cite[Theorem 9]{GP} that an endomorphism $\Phi$ of $\cS$ is
an automorphism of the Lie algebra $\cS$ if and only if it has the form
\[\Phi=\overline{\phi}\circ{\cal U} _{K}\circ e^{\overline{\Omega}}.\]
Here $\phi\in \op{Diff}(M)$ and if $\cS$ is interpreted as the
algebra $\op{Pol}(T^*M)$ of polynomial functions on $T^*M$, the
automorphism $\overline{\phi}$ is implemented by the phase lift
$\phi^*$ of $\phi$ to the cotangent bundle $T^*M$, a
symplectomorphism of $T^*M$. If, on the other hand, $\cS$ is
viewed as the algebra $\Gamma(\cS TM)$ of symmetric contravariant
tensor fields on $M$, the automorphism $\overline{\phi}$ is the
standard action of $\phi$ on such tensor fields.

Further, $K\in\R^*$, $\Omega$ is a closed 1-form on $M$ and the
automorphism ${\cal U}_K\in \op{Aut}\cS$, ${\cal U}_K:\cS_i\ni
S\raa K^{i-1}S\in\cS_i$, is, for $K>0$, the exponential of the
derivation $\ln{K}\,\op{Deg}$, whereas the automorphism
$e^{\overline{\Omega}}$ induced by the lowering derivation
$\overline{\Omega}$, i.e. the action of the vertical vector field
${\Omega}^v$, is the composition with the translation $\cT_{\zW}$
by $\zW$ in $T^*M$. Note that since the homothety $h_K$ of $T^*M$
by $K$ acts on homogeneous polynomials of degree $i$  by
multiplication by $K^i$, the automorphism $\cU_K$ can be written
also in the form $\cU_K(S)=K^{-1}\,S\circ h_K$. Hence, every
one-parameter group of automorphisms of the Lie algebra $\cS$ has
the form
$$\zF_{\phi_t,K_t,\zW_t}(S)=K_t^{-1}\, S\circ\cT_{\zW_t}\circ
h_{K_t}\circ(\phi_t^{-1})^*.
$$
It is easy to prove the following commutation relations.
\begin{prop} \beas h_K\circ\phi^*&=&\phi^*\circ h_K\ ,\\
\cT_\zW\circ h_K&=&h_K\circ\cT_{K^{-1}\,\zW}\ ,\\
(\phi^{-1})^*\circ\cT_{\zW}&=&\cT_{\phi^*\zW}\circ(\phi^{-1})^*.\eeas
\end{prop}
These relations together with the one-parameter group property
yield
\beas\phi_{t+s}&=&\phi_t\circ\phi_s\ ,\\
K_{t+s}&=&K_t\cdot K_s\ ,\\
\zW_{t+s}&=&\zW_t+K_t\cdot\phi^*_t\zW_s \ , \eeas with the initial
conditions $\phi_0=\op{id}_M$, $K_0=1$, $\zW_0=0$. The obvious
unique solutions are $\phi_t=\op{Exp}(tY)$ for a certain complete
vector field $Y$, $K_t=e^{\zk t}$ for a certain $\zk\in\R$, and
$$\zW_t=\int^t_0e^{\zk s}\, (\op{Exp}(sY))^*\zw\,ds
$$
for a certain closed 1-form $\zw$ on $M$.\\

Let us systematically characterize derivations by the unique
triplet with first member vanishing on the $0$-section. As
well-known there is a Lie algebra isomorphism between
${\cS}_1(M)=\op{Pol}^1(T^*M)$ and ${\cX}(M).$ We denote a
homogeneous first order polynomial $P$ on $T^*M$ by $P_*$ when it
is viewed as vector field of $M$. Note that the hamiltonian vector
field $X_P$ of $P$ is nothing but the phase lift $(P_*)^*$ of
$P_*$. We then have the following theorem.
\begin{theo} A derivation
\begin{equation}C_{P,\zk,\zw}(S)=\{P,S\}+\kappa\,\op{Deg}(S)+\zw^v(S)
\end{equation}
of the Lie algebra $\cS(M)$ of all infinitely differentiable
functions of $T^*M$ that are polynomial along the fibers, where
$P$ is vanishing on the $0$-section, is integrable if and only if
the polynomial function $P$ belongs to $\cS_1(M)$ and is complete,
i.e. the hamiltonian vector field $X_P$ of $P$ is complete, i.e.
the basis vector field $P_*$ is complete. In this case the
one-parameter group of automorphisms $\zF_t$ generated by
$C_{P,\zk,\zw}$ reads
$$\zF_t(S)=e^{-\zk t}\, S\circ\cT_{\int^t_0e^{\zk s}\,
(\op{Exp}(sP_*))^*\zw\,ds}\circ h_{e^{\zk t}}\circ
\op{Exp}(-tX_P).
$$
\end{theo}

\subsection{The case $\cD(M)$}

Let us eventually recall that Theorem $10$ of \cite{GP} asserts
that automorphisms of the Lie algebra $\cD$ have the form
\be\label{C}\Phi=\phi_*\circ\cC^{a}\circ e^{\overline{\Omega}},
\ee where $e^{\overline{\Omega}}$ ($\zW\in\zW^1(M)\cap
\op{ker}\,d$) is the formerly mentioned automorphism of $\cD$ and
where $\phi_*$ $(\phi\in \op{Diff}(M))$ is the automorphism of
$\cD$ defined by $\phi_*(D)=\phi\circ D\circ\phi^{-1},$ $\phi(f)$
being of course $f\circ\phi^{-1}$ ($D\in\cD,f\in\cA$). Moreover,
superscript $a$ is $0$ or $1$, so that $\cC^{a}$ is
$\cC^0=\op{id}$ or $\cC^1=\cC$, $\cC$ being the opposite of the
conjugation operator $*$. Remember that for an oriented manifold
$M$ with volume form $\eta$, the conjugate $D^*\in\cD$ of a
differential operator $D\in\cD$ is defined by
\[\int_MD(f)\cdot g\mid\eta\mid=\int_Mf\cdot D^*(g)\mid\eta\mid,\]
for any compactly supported $f,g\in\cA$. Since
$(D\circ\De)^*=\De^*\circ D^*$ $(D,\De\in\cD)$, the operator
$\cC:=-*$ verifies $\cC(D\circ\De)=-\cC(\De)\circ\cC(D)$ and is
thus an automorphism of $\cD$. Formal calculus allows to show that
this automorphism exists for any manifold, orientable or not.
Clearly, the automorphism $\cC$ is not implemented by a derivation
and (\ref{C}) belongs to the connected component of identity only
if $a=0$. Thus we can consider one-parameter groups of
automorphisms of the form
$$\zF_{\phi_t,\zW_t}=(\phi_t)_*\circ e^{\overline{\Omega_t}}.
$$
It is easy to prove that
$$e^{\overline{\Omega}}\circ\phi_*=\phi_*\circ
e^{\overline{\phi^*\Omega}},
$$
so that the one-parameter group property yields
\beas \phi_{t+s}&=&\phi_t\circ\phi_s\ ,\\
\zW_{t+s}&=&\zW_t+\phi_t^*\zW_s\ , \eeas with initial conditions
$\phi_0=\op{id}_M$, $\zW_0=0$. The obvious general solutions are
$\phi_t=\op{Exp}(tY)$ for a complete vector field $Y$, and
$\zW_t=\int_0^t(\op{Exp}(sY))^*\zw\,ds$ for a certain closed
1-form $\zw$. Thus we get the following.

\begin{theo} A derivation
\begin{equation}C_{P,\zw}(D)=[P,D]+\overline{\zw}(D)
\end{equation}
of the Lie algebra $\cD(M)$ of all differential operators is
integrable if and only if $P\in{\cX}(M)$ and $P$ is complete. In
this case the one-parameter group of automorphisms $\zF_t$
generated by $C_{P,\zw}$ reads
$$\zF_t=(\op{Exp}(tP))_*\circ\,e^{\overline{\int_0^t(\op{Exp}(sP))^*\zw\,ds}}.
$$
\end{theo}

\medskip\noindent
{\bf Remark:} The results of this section describing commutations
rules for automorphisms easily imply that $\op{Aut}\cP$ is an
infinite-dimensional regular Lie group in the sense of A.~Kriegl
and P.~Michor (see \cite{KM1} or \cite[Ch. 8]{KM2}). The
integrable derivations (in fact, those with compact supports) form
the Lie algebra of $\op{Aut}\cP$.

\noindent Janusz GRABOWSKI\\Polish Academy of Sciences\\Institute
of Mathematics\\\'Sniadeckich 8\\P.O. Box 21\\00-956 Warsaw,
Poland\\Email: jagrab@impan.gov.pl\\\\
\noindent Norbert PONCIN\\University of Luxembourg\\Mathematics Laboratory\\avenue
de la Fa\"{\i}encerie, 162A\\
L-1511 Luxembourg City, Grand-Duchy of Luxembourg\\Email:
norbert.poncin@uni.lu
\end{document}